\documentclass{article}
 \usepackage{amsmath}

\newtheorem{thm}{Theorem}[section]

\newtheorem{prop}[thm]{Proposition}

\newcommand{\be}{\begin{equation}}
\newcommand{\ee}{\end{equation}}
\newcommand{\ben}{\begin{enumerate}}
\newcommand{\een}{\end{enumerate}}
\newcommand{\beq}{\begin{eqnarray}}
\newcommand{\eeq}{\end{eqnarray}}
\newcommand{\beqn}{\begin{eqnarray*}}
\newcommand{\eeqn}{\end{eqnarray*}}

\newcommand{\pa}{\partial}

\newcommand{\qed}{\hspace*{\fill}Q.E.D.}  

\begin{document}
\title{On a Class of Singular Projectively Flat Finsler Metrics with Constant Flag Curvature}
\author{Guojun Yang}
\date{}
\maketitle
\begin{abstract}
Singular Finsler metrics, such as Kropina metrics and $m$-Kropina
metrics, have a lot of applications in the real world. In this
paper, we classify a class of  singular $(\alpha,\beta)$-metrics
which are locally projectively flat with constant flag curvature
in dimension $n=  2$ and $n \ge 3$ respectively. Further, we
determine the local structure of  $m$-Kropina metrics and
particularly Kropina metrics  which are projectively flat with
constant flag curvature and prove that such metrics must be
locally Minkowskian but  are not necessarily flat-parallel.

\bigskip


\end{abstract}

\section{Introduction}

In Finsler geometry, the flag curvature is a natural extension of
the sectional curvature in Riemannian geometry. The flag curvature
of a Finsler metric on a manifold $M$ is a scalar function ${\bf
K}= {\bf K}(x, y, P)$ of a tangent plane $P\subset T_xM$  and a
non-zero vector $y\in P $. It is said to be of constant flag
curvature if ${\bf K}$ is a constant, and isotropic flag curvature
if ${\bf K}={\bf K}(x)$. The Schur Theorem shows that ${\bf
K}=constant$ if ${\bf K}={\bf K}(x)$ and the dimension $n\ge 3$.
The Beltrami Theorem says that a Riemannian metric is of constant
sectional curvature if and only if it is locally projectively
flat, that is, geodesics are straight lines locally. However,
locally projectively flat Finsler metrics are not necessarily of
constant flag curvature or isotropic flag curvature (generally are
of scalar flag curvature, namely, ${\bf K}={\bf K}(x, y)$
independent of $P$). On the other hand, it is easy to prove that
any two-dimensional projectively flat Finsler metric with
isotropic flag curvature is of constant flag curvature.

 An {\it $(\alpha,\beta)$-metric} is defined by a
Riemannian metric
 $\alpha=\sqrt{a_{ij}(x)y^iy^j}$ and a $1$-form $\beta=b_i(x)y^i$ on a manifold
 $M$, which can be expressed in the following form:
 $$F=\alpha \phi(s),\ \ s=\beta/\alpha,$$
where $\phi(s)$ is a $C^{\infty}$ function on $(-b_o,b_o)$. It is
known that $F$ is a regular Finsler metric for any $(\alpha,
\beta)$ with $ \|\beta\|_{\alpha} < b_o$ if and only if
 \be\label{j1}
 \phi(s)>0,\quad \phi(s)-s\phi'(s)+(\rho^2-s^2)\phi''(s)>0,
   \quad (|s| \leq \rho <b_o),
 \ee
 where $b_o$ is a constant (\cite{Shen2}). In this paper, we do
  not assume  the regular condition, and we will study a class of
  singular $(\alpha,\beta)$-metrics satisfying (\ref{j02}) below.  Singular Finsler metrics have a lot
of applications in the real world (\cite{AIM} \cite{AHM} ). Z.
Shen also introduces singular Finsler metrics in \cite{Shen3}.

A {\it Randers metric} is in the form $F= \alpha+\beta$, which is
a special $(\alpha,\beta)$-metric. It is proved in \cite{Shen0}
that a Randers metric
 $F=\alpha+\beta$ is projectively flat with constant flag
 curvature if  and only if $F$ is locally Minkowskian (equivalently {\it flat-parallel}, that is,
 $\alpha$ is flat and $\beta$ is parallel with respect to $\alpha$) or after scaling,
 $\alpha$ and $\beta$ can be locally
  written as
  \be\label{j2}
  \alpha=\frac{\sqrt{(1-|x|^2)|y|^2+\langle x,y\rangle^2}}{1-|x|^2}, \ \
  \beta=\pm\Big\{\frac{\langle x,y\rangle}{1-|x|^2}
   +\frac{\langle a,y\rangle}{1-|x|^2}\Big\},
  \ee
where $a\in R^n$ is a constant vector. The Randers metric
$F=\alpha+\beta$ defined by (\ref{j2}) is of constant flag
curvature $K=-1/4$.

It is proved in \cite{SY} that an $(\alpha,\beta)$-metric in the
form $F=(\alpha+\beta)^2/\alpha$ is projectively flat with
constant flag curvature if and only if $F$ is locally Minkowskian
(equivalently flat-parallel) or after scaling, $\alpha$ and
$\beta$ can be locally written as,
  \be\label{j3}
 \alpha=\lambda\frac{\sqrt{(1-|x|^2)|y|^2+\langle x,y\rangle^2}}{1-|x|^2},\ \
  \beta=\pm\lambda\Big\{\frac{\langle x,y\rangle}{1-|x|^2}
   +\frac{\langle a,y\rangle}{1-|x|^2}\Big\},
  \ee
where
$$\lambda:=\frac{(1+\langle a,x\rangle)^2}{1-|x|^2}.$$
The  metric $F=(\alpha+\beta)^2/\alpha$ defined in (\ref{j3}) is
of zero flag curvature (see \cite{MSY}).

In \cite{Shen3}, Z. Shen gives the Taylor expansions for
$x$-analytic projectively flat metrics $F=F(x,y)$ of constant flag
curvature $K$. For some suitable choices of $K$,
$\psi(y)(=F|_{x=0})$ and $\varphi(y)(=F_{x^k}y^k/(2F)|_{x=0})$,
one can easily get the projectively flat Finsler metrics of
constant flag curvature ${\bf K}=K$ in (\ref{j2}) and (\ref{j3}).

In \cite{LS1}, the authors classify projectively flat
$(\alpha,\beta)$-metrics of constant flag curvature in dimensions
$n\ge 3$ and $\phi(0)=1$. They show  that such $(\alpha,
\beta)$-metrics must be flat-parallel, or after a suitable
scaling, isometric to the metrics in (\ref{j2}) or (\ref{j3}).
When $n=2$ and $\phi(0)=1$, the present author shows the
essentially same conclusions (see \cite{Y2}).

For an $(\alpha,\beta)$-metric $F=\alpha \phi(\beta/\alpha)$, if
  $\phi(0)> 0$, then generally $F$ is regular. However, if $\phi(0)= 0$  or $\phi(0)$ is not defined,
  then  $\phi$ does not satisfy (\ref{j1}) and in this
  case $F$ is singular.
  In this paper,  we assume
  $\phi(s)$ is in the following form
   \be\label{j02}
 \phi(s):=cs+s^m\varphi(s),
   \ee
where $c,m$ are constant with $m\ne 0,1$ and $\varphi(s)$ is a
$C^{\infty}$ function on a neighborhood of $s=0$ with
$\varphi(0)=1$, and further for convenience we put $c=0$ if $m$ is
a negative integer. If $m=0$, we have $\phi(0)=1$ and this case
appears in a lot of literatures. When $m\ge 2$ is an integer,
(\ref{j02})
  is equivalent to the following condition
  $$
  \phi(0)=0,\quad
  \phi^{(k)}(0)=0 \ \ (2\le k\le m-1), \quad \phi^{(m)}(0)=m!.
  $$
  Another interesting case is $c=0$ and $\varphi(s)\equiv 1$ in
  (\ref{j02}), and in this case, $F=\alpha\phi(s)$ is called an
  $m$-Kropina metric, and in particular a Kropina metric when $m=-1$.

In \cite{Y3} \cite{Y4},   the present author classifies the
$(\alpha,\beta)$-metric $F=\alpha\phi(\beta/\alpha)$ which is
Douglasian and locally projectively flat respectively, where
$\phi(s)$ satisfies (\ref{j02}). In this paper we prove the
following theorem:

\begin{thm}\label{th1}
  Let $F=\alpha \phi(s)$, $s=\beta/\alpha$, be an
  $(\alpha,\beta)$-metric on an open subset $U$ of the $n$-dimensional
  Euclidean space $R^n$, where $\phi(s)$ satisfies (\ref{j02}). Then $F$ is projectively
  flat with constant flag curvature $K$ if and only if one
  of the following cases holds
  \ben
    \item[{\rm (i)}] ($n\ge 2$) $F$ is flat-parallel, namely, $\alpha$ is  flat and $\beta$ is
    parallel  with respect to $\alpha$.

     \item[{\rm (ii)}] ($n\ge 2$) $F$ is given by
    \be\label{Fii}
      F=\beta^m(\alpha^2+k\beta^2)^{\frac{1-m}{2}},
    \ee
    which is  projectively flat with  $K=0$, where
     $k$ is a constant. In this case, $F$ is locally Minkowskian,
     but generally are not flat-parallel.

     \item[{\rm (iii)}] ($n=2$) $F$ is given by
     \be
         F=\frac{2b}{1-kb^2}\Big\{b\sqrt{\alpha^2-k\beta^2}
         -\sqrt{b^2\alpha^2-\beta^2}\Big\},\label{Fiii}
     \ee
     which is projectively flat with $K<0$, where $b=||\beta||_{\alpha}$
     and $k$ are two constants with $k\ne 1/b^2$.

     \item[{\rm (iv)}] ($n=2$) $F$ is given by
    \be
      F=\frac{4b^2}{(1-kb^2)^2}\frac{\big(b\sqrt{\alpha^2-k\beta^2}
          -\sqrt{b^2\alpha^2-\beta^2}\big)^2}{\sqrt{\alpha^2-k\beta^2}},\label{Fiv}
   \ee
   which  is projectively flat with $K=0$, where $b=||\beta||_{\alpha}$
     and $k$ are two constants with $k\ne 1/b^2$.
  \een
\end{thm}

The metric  $F$ in Theorem \ref{th1}(i) is a special locally
Minkowskian metric.  Since $F$ is never Riemannian, an easy proof
shows that $\alpha$ is flat if $\beta$ is parallel with respect to
$\alpha$.

Note that in dimension $n=2$, the term
$\sqrt{b^2\alpha^2-\beta^2}$ is actually a $1$-form. Therefore,
the metric in (\ref{Fiii}) is essentially a Randers metric
$\widetilde{F}=\widetilde{\alpha}+\widetilde{\beta}$ when
$k<1/b^2$, where
 $$\widetilde{\alpha}:=\frac{2b^2}{1-kb^2}\sqrt{\alpha^2-k\beta^2},
\ \
\widetilde{\beta}:=-\frac{2b}{1-kb^2}\sqrt{b^2\alpha^2-\beta^2}.$$
The metric in (\ref{Fiv}) is essentially the type
$\widetilde{F}=(\widetilde{\alpha}+\widetilde{\beta})^2/\widetilde{\alpha}$
when $k<1/b^2$, where
$$\widetilde{\alpha}:=\frac{4b^4}{(1-kb^2)^2}\sqrt{\alpha^2-k\beta^2},
\ \
\widetilde{\beta}:=-\frac{4b^3}{(1-kb^2)^2}\sqrt{b^2\alpha^2-\beta^2}.$$
Thus according to \cite{MSY}, \cite{Shen0} and \cite{SY}, the
local structures of the $(\alpha,\beta)$-metrics in (\ref{Fiii})
and (\ref{Fiv}) can be determined.

\bigskip

 As seen above,
$F=\alpha+\beta$ or $F=(\alpha+\beta)^2/\alpha^2$ is locally
Minkowskian if and only if $F$ is flat-parallel. However, the
metric in (\ref{Fii}) is not necessarily the case. When $k>
-1/b^2$, the  metric in (\ref{Fii}) is essentially an $m$-Kropina
metric,
 $$
 F=\bar{\beta}^m\bar{\alpha}^{1-m}, \ \ \
 (\bar{\alpha}:=\sqrt{\alpha^2+k\beta^2},\ \bar{\beta}:=\beta).
 $$
Now for an $m$-Kropina metric (we may put $k=0$ in (\ref{Fii})),
we can determine its local structure as follows.

\begin{thm}\label{th01}
 Let $F=\beta^m\alpha^{1-m}$ be an $m$-Kropina metric, where $m\ne
 0,1$. Suppose $F$ is locally projectively flat with vanishing
 flag curvature. Then $F$ can be written as
 $F=\widetilde{\alpha}^{1-m}\widetilde{\beta}^m$,
  where $\widetilde{\alpha}$ is flat and $\widetilde{\beta}$
  is parallel with respect to $\widetilde{\alpha}$, and thus $\widetilde{\alpha}$ and
  $\widetilde{\beta}$ can be locally written as
 \be\label{ycw1}
\widetilde{\alpha}=|y|,\ \ \widetilde{\beta}=y^1.
 \ee
Further $\alpha,\beta$ are related
    with $\widetilde{\alpha},\widetilde{\beta}$ by
    \be\label{ycw2}
  \alpha=\eta^{\frac{m}{m-1}}\widetilde{\alpha}, \ \ \ \beta=\eta
  \widetilde{\beta},
    \ee
where $\eta=\eta(x)>0$ is a scalar function. Obviously $F$ is
locally Minkowskian, but generally  not flat-parallel.
\end{thm}

 In \cite{NLB}, the authors claim that for a projectively flat
 Kropina metric $F=\alpha^2/\beta$ with vanishing flag curvature,
 $\alpha$ must be flat and $\beta$ must be closed. However, Theorem \ref{th01} shows a different conclusion. By
 (\ref{ycw2}), generally $\alpha$ is not flat and $\beta$ is not
 closed since $\eta$ can be arbitrary.

In Theorem \ref{th01}, if $m\ne -1$, we can obtain the same
conclusion under weaker conditions---$F$ is only assumed to be
Douglasian if $n=2$ (see \cite{Y4}), or $F$ is only assumed to be
locally projectively flat if $n\ge 3$ (see \cite{Y3}), or $F$ is
only assumed to be of constant/scalar flag curvature for $n\ge 3$
(\cite{Y5}).

  The general characterization for the
metric $F$ in (\ref{Fii}) which is locally projectively flat with
vanishing flag curvature is given by the equations
(\ref{y16})--(\ref{w03}) below with $P=0$, where we should note
that if we put $m=-1$ in (\ref{y16})--(\ref{w03}), then we get
(\ref{ygjcw}), (\ref{w001}), (\ref{cr5}) and (\ref{cr7}) with
$\mu=-2b^2\tau$. But it seems difficult to obtain their local
solutions.

\section{Preliminaries}

  In local coordinates, the geodesics of a Finsler metric
  $F=F(x,y)$ are characterized by
   $$\frac{d^2 x^i}{d t^2}+2G^i(x,\frac{d x^i}{d t})=0,$$
where
 \beq \label{G1}
 G^i:=\frac{1}{4}g^{il}\big \{[F^2]_{x^ky^l}y^k-[F^2]_{x^l}\big \}.
 \eeq
 The local functions $G^i$ are called the spray coefficients of
 $F$.
A Finsler metric $F$ is said to be projectively flat in $U$, if
 there is a local coordinate system $(U,x^i)$ such that
  $G^i=Py^i$,
  where $P=P(x,y)$ is called the  projective factor. In this case,
  the scalar
  flag curvature $K$ is given by
  \be\label{y4}
  K=\frac{P^2-P_{x^k}y^k}{F^2}.
  \ee

  Consider an $(\alpha,\beta)$-metric $F=\alpha \phi(s),
s=\beta/\alpha$. Let $\nabla \beta = b_{i|j} y^i dx^j$ denote the
covariant derivatives of $\beta$ with respect to $\alpha$ and
define
 $$r_{ij}:=\frac{1}{2}(b_{i|j}+b_{j|i}),\ \ s_{ij}:=\frac{1}{2}(b_{i|j}-b_{j|i}),\ \
  s_j:=b^is_{ij},\ \ s^i:=a^{ik}s_k,$$
 where $b^i:=a^{ij}b_j$ and $(a^{ij})$ is the inverse of
 $(a_{ij})$. Then  by (\ref{G1}), the spray coefficients
$G^i$ of $F$
 are given by (\cite{CS}  \cite{LSS} \cite{LS2} \cite{Ma}  \cite{Shen1} \cite{Shen2}):
  \be\label{y005}
  G^i=G^i_{\alpha}+\alpha Q s^i_0+\alpha^{-1}\Theta (-2\alpha Q
  s_0+r_{00})y^i+\Psi (-2\alpha Q s_0+r_{00})b^i,
  \ee
where $s^i_j=a^{ik}s_{kj}, s^i_0=s^i_ky^k,
s_i=b^ks_{ki},s_0=s_iy^i$, and
 $$
  Q:=\frac{\phi'}{\phi-s\phi'},\ \
  \Theta:=\frac{Q-sQ'}{2\Delta},\ \
  \Psi:=\frac{Q'}{2\Delta},\ \ \Delta:=1+sQ+(b^2-s^2)Q'.
 $$
By (\ref{y005}), it is easy to see that if $\alpha$ is
projectively flat and $\beta$ is parallel with respect to $\alpha$
($r_{ij}=0,s_{ij}=0$), then $F=\alpha\phi(\beta/\alpha)$ is
projectively flat.

\bigskip

In this paper, our proof is based on the following theorem.

\begin{thm} {\rm (\cite{Y3} \cite{Y4})}\label{th4}
  Let $F=\alpha \phi(s)$, $s=\beta/\alpha$, be an
  $(\alpha,\beta)$-metric on an open subset $U\subset R^n$.
   Suppose that $\beta$ is not parallel with respect to
     $\alpha$ and $\phi$ satisfies (\ref{j02}). Let $G^i_{\alpha}$
     be the spray coefficients of $\alpha$.
  Then $F$ is  projectively flat in $U$ with
 $G^i=P(x,y)y^i$ if and only if one of the following cases holds:

 \ben
 \item[{\rm (i)}] ($n\ge 2$) For a 1-form $\rho=\rho_i(x)y^i$, $\phi(s)$,  $\beta$ and $G^i_{\alpha}$ satisfy
 \beq
\phi(s)&=&ks+\frac{1}{s}, \ \ \ \ s_{ij}=\frac{b_is_j-b_js_i}{b^2},\label{ygjcw}\\
  G^i_{\alpha}&=&\rho
y^i-\frac{r_{00}}{2b^2}b^i-\frac{\alpha^2-k\beta^2}{2b^2}s^i,\label{w001}
 \eeq
where $k$ is a constant. In this case, the projective factor $P$
is given by
 \be\label{w0001}
 P=\rho-\frac{1}{b^2(\alpha^2+c\beta^2)}\Big\{(\alpha^2-c\beta^2)s_0+r_{00}\beta\Big\}.
 \ee

  \item[{\rm (ii)}] ($n\ge 2$) For a 1-form $\rho=\rho_i(x)y^i$ and a scalar $\tau=\tau(x)$, $\phi(s)$, $\beta$ and $G^i_{\alpha} $satisfy
   \beq
    \phi(s)&=&a_1s+s^m(1+k s^2)^{\frac{1-m}{2}},\label{y5}\\
    b_{i|j}&=&2\tau \big\{mb^2a_{ij}-(m+1+kb^2)b_ib_j\big\},\label{cw6}\\
   G^i_{\alpha}& =&\rho
   y^i-\tau(m\alpha^2-k\beta^2)b^i,\label{cw7}
   \eeq
  where $a_1$ and $k$ are constant. In this case, the projective factor $P$ is given by

  \be\label{w01}
  P=\rho +\tau \alpha
  \Big\{s(-m+ks^2)-s^2(1+ks^2)\frac{\phi'}{\phi}\Big\}
  \ee

\item[{\rm (iii)}] ($n\ge2$) For a 1-form $\rho=\rho_i(x)y^i$ and
a scalar $\tau=\tau(x)$, $\phi(s)$,  $\beta$ and $G^i_{\alpha}$
satisfy
 \beq
\phi(s)&=&s^m(1+k s^2)^{\frac{1-m}{2}},\ \ \ \ s_{ij}=\frac{b_is_j-b_js_i}{b^2},\label{y16}\\
 r_{ij}&=&2\tau \big\{mb^2a_{ij}-(m+1+kb^2)b_ib_j\big\}
   -\frac{m+1+2kb^2}{(m-1)b^2}(b_is_j+b_js_i),\label{w06}\\
   G^i_{\alpha}&=&\rho
   y^i+\Big\{\frac{2k\beta s_0}{(m-1)b^2}-\tau(m\alpha^2-k\beta^2)\Big\}b^i
   -\frac{m\alpha^2+k\beta^2}{(m-1)b^2}s^i,\label{ycr06}
 \eeq
 where  $k$ is a constant. In this case, the projective factor $P$ is given by
  \be\label{w03}
  P=\rho-2m\tau\beta-\frac{2m}{(m-1)b^2}s_0.
  \ee

\item[{\rm (iv)}] ($n= 2$) For a 1-form $\rho=\rho_i(x)y^i$ and a
scalar $\tau=\tau(x)$, $\phi(s)$,  $\beta$ and $G^i_{\alpha}$
satisfy
 \beq
\phi(s)&=&mb^2\sqrt{b^2-s^2}\int_0^s\frac{1}{(b^2-t^2)^{3/2}}
    \Big(\frac{t}{\sqrt{1-k t^2}}\Big)^{m-1}dt,\label{j4}\\
    r_{ij}&=&-\frac{1}{b^2}(b_is_j+b_js_i),\label{j5}\\
    G^i_{\alpha}&=&\rho
   y^i-\frac{(m-2)\alpha^2+k\beta^2}{(m-1)b^2}s^i,\label{j6}
 \eeq
 where  $k$ is a constant. In this case, the projective factor $P$ is given by
  \be\label{w04}
  P=\rho+\frac{1}{(m-1)b^2}\Big\{s(ks^2-1)\frac{\phi'}{\phi}-ks^2-m+2\Big\}s_0.
  \ee

   \item[{\rm (v)}]  ($n= 2$)  For a 1-form $\rho=\rho_i(x)y^i$ and a scalar $\tau=\tau(x)$, $\phi(s)$ and $\beta$ satisfy
\beq
    \phi(s)&=&k_1s+\frac{2k_2}{s}+\frac{1}{s^3},\label{yg5} \\
   r_{ij}&=&-2\tau
   \big\{3b^2a_{ij}+(k_2b^2-2)b_ib_j\big\}+\frac{(3k_1+k_2^2)b^4-4}{8b^2(1+k_2b^2)}(b_is_j+b_js_i),\label{yg6}\\
   G^i_{\alpha}&=&\rho
y^i+\tau(3\alpha^2+k_2\beta^2)b^i+\Big\{\frac{k_1-k_2^2}{8(1+k_2b^2)}(3b^2\alpha^2-\beta^2)\nonumber\\
&&+(\frac{k_2}{2}-\frac{3}{4b^2})\alpha^2-\frac{k_2}{b^2}\beta^2\Big\}s^i,\label{cw002}
   \eeq
 where  $k_1,k_2$ are
constant satisfying $1+k_2b^2\ne 0$.  In this case, the projective
factor $P$ is given by
 \beq\label{cw0002}
 P&=&\rho+2\tau\beta\Big\{3-\frac{2c\beta^4}{\alpha^4+c\beta^4+k_2\beta^2(2\alpha^2+k_2\beta^2)}\Big\}
 +(\frac{k_2b^2-3}{2b^2}+T)s_0,\\
 T:&=&c\frac{4\beta^2(2\beta^2-b^2\alpha^2)+3b^4(\alpha^4+c\beta^4)+k_2b^2\beta^2(6b^2\alpha^2+4\beta^2+3k_2b^2\beta^2)}
 {8b^2(1+k_2b^2)\big[\alpha^4+c\beta^4+k_2\beta^2(2\alpha^2+k_2\beta^2)\big]}.\nonumber\\
 c:&=&k_1-k_2^2.\nonumber
 \eeq
  \een

\end{thm}

\section{The first class in Theorem \ref{th4}}\label{first}

In this section, we study the property of the
$(\alpha,\beta)$-metric determined by (\ref{ygjcw}) and
(\ref{w001}) in Theorem \ref{th4}.

\begin{prop}
 Let $F=k\beta+\alpha^2/\beta$, where $k$ is a constant, be an $n$-dimensional
 ($\alpha,\beta$)-metric which is  projectively flat with
 constant flag curvature $K$.   Then we have  $K=0$, and then $F$
 is locally Minkowskian.
\end{prop}

{\bf Proof.} We only need to assume that $\beta$ is not parallel
with respect to $\alpha$. The projective factor $P$ is given by
(\ref{w0001}). By (\ref{w001}) we get
 \be\label{cw23}
 \alpha_0:=\alpha_{x^k}y^k=a_{ir}y^r\frac{2}{\alpha}G^i_{\alpha}=2\rho\alpha+\frac{(cs^2-1)s_0\alpha-sr_{00}}{b^2}.
 \ee
Now it follows from (\ref{w001}) and (\ref{cw23}) that there holds
 \be
 s_{x^k}y^k=-\frac{s\alpha_0}{\alpha}+\frac{r_{00}+2b_iG^i_{\alpha}}{\alpha}
           =-\frac{s\big[(cs^2-1)s_0\alpha-sr_{00}\big]}{b^2\alpha}.\label{cw24}
 \ee

Then plug $\phi(s)=ks+1/s$,  (\ref{w0001}), (\ref{cw23}) and
(\ref{cw24}) into (\ref{y4}) and thus (\ref{y4}) can be written in
the following form
 \be\label{cr1}
 A_1\beta^2+Kb^4\alpha^8=0,
 \ee
where $A_1$ is a homogeneous polynomial in $(y^i)$ of degree six.
Clearly by (\ref{cr1}) we get $K=0$. By $K=0$, (\ref{cr1}) has the
following equivalent form
 \be\label{cr4}
 A_2(\alpha^2+k\beta^2)+3\beta^2(r_{00}-2k\beta s_0)^2=0,
 \ee
where  $A_2$ is a polynomial in $(y^i)$. So (\ref{cr4}) implies
 \be\label{cr5}
 r_{00}=2k\beta s_0+\mu(\alpha^2+k\beta^2),
 \ee
where $\mu=\mu(x)$ is a scalar function. Now plug (\ref{cr5}) into
(\ref{w0001}) and then we have
 \be\label{cr6}
 P=\rho-\frac{\mu\beta+s_0}{b^2}.
 \ee
So $P$ is a 1-form and thus $G^i=Py^i$ are quadratic, which shows
$F$ is Berwaldian. Plus $K=0$, $F$ is locally Minkowskian. So we
have $P=0$ and thus again by (\ref{cr6}) we have
 \be\label{cr7}
 \rho=\frac{\mu\beta+s_0}{b^2}.
\ee

\section{The second class in Theorem \ref{th4}}\label{sec3}
In this section, we study the property of the
$(\alpha,\beta)$-metric determined by (\ref{y5}), (\ref{cw6}) and
(\ref{cw7}) in Theorem \ref{th4}.

\begin{prop}
 Let $F$ be an $n$-dimensional ($\alpha,\beta$)-metric given by (\ref{y5}) which is  projectively flat with
 constant flag curvature $K$. Suppose $\beta$ is not parallel with respect to $\alpha$. Then we have $a_1=0$, $K=0$.
\end{prop}

{\bf Proof.}
 The projective factor $P$ is given by (\ref{w01}). By (\ref{cw6}) we
get
 \be\label{y22}
 r_{00}=2\tau\big(mb^2-(m+1+kb^2)s^2\big)\alpha^2.
 \ee
By (\ref{cw7}) we get
 \be\label{y23}
 \alpha_0:=\alpha_{x^k}y^k=a_{ir}y^r\frac{2}{\alpha}G^i_{\alpha}=2\tau s(ks^2-m)\alpha^2+2\rho\alpha.
 \ee
Now it follows from (\ref{cw7}), (\ref{y22}) and (\ref{y23}) that
there holds
 \be
 s_{x^k}y^k=-\frac{s\alpha_0}{\alpha}+\frac{r_{00}+2b_iG^i_{\alpha}}{\alpha}
           =-2\tau\alpha s^2(1+ks^2).\label{y24}
 \ee
The function $\phi(s)$ in (\ref{y5}) satisfies
 \be\label{y024}
 \phi''=\frac{(ks^2-m)(\phi-s\phi')}{s^2(1+ks^2)}.
 \ee
Then plug (\ref{w01}), (\ref{y23})--(\ref{y024}) into (\ref{y4})
and we obtain

 \beq
 K\alpha^2\phi^2\hspace{-0.6cm}&&=\rho^2-\rho_{00}+(m+1)\big(\tau_0-(m+1)s\tau^2\alpha\big)s\alpha+
 3\frac{(\phi-s\phi')^2}{\phi^2}s^2(1+ks^2)^2\tau^2\alpha^2\nonumber\\
 &&
 +\frac{\phi-s\phi'}{\phi}s(1+ks^2)\big[2(m-1)\tau^2s\alpha-\tau_0\big]\alpha,\label{y25}
 \eeq
where
 \be\label{y26}
 \tau_i:=\frac{\pa \tau}{\pa x^i}, \ \ \tau_0:=\tau_iy^i,\ \ \rho_{ij}:=\frac{1}{2}(\frac{\pa \rho_i}{\pa x^j}+\frac{\pa
  \rho_j}{\pa x^i}),\ \ \rho_{00}:=\rho_{ij}y^iy^j.
 \ee

To deal with (\ref{y25}), we choose a special coordinate system at
a point as that in \cite{Shen1}. At a fixed point $x_o$, make a
change of coordinates: $(s,y^A)\mapsto (y^1,y^A)$ by
 $$y^1=\frac{s}{\sqrt{b^2-s^2}}\bar{\alpha},\ \ y^A=y^A,$$
where $\bar{\alpha}=\sqrt{\sum_{A=2}^n (y^A)^2}$. Then
$$\alpha=\frac{b}{\sqrt{b^2-s^2}}\bar{\alpha},\ \
\beta=\frac{bs}{\sqrt{b^2-s^2}}\bar{\alpha}.$$
 We get
 \be\label{y27}
  \rho=\frac{s\rho_1}{\sqrt{b^2-s^2}}\bar{\alpha}+\bar{\rho}_0,
  \ \ \rho_{00}=\frac{\rho_{11}s^2}{b^2-s^2}\bar{\alpha}^2
  +\frac{2s\bar{\rho}_{10}}{\sqrt{b^2-s^2}}\bar{\alpha}+\bar{\rho}_{00},\ \
  \tau_0=\frac{s\tau_1}{\sqrt{b^2-s^2}}\bar{\alpha}+\bar{\tau}_0,
  \ee
 where $\bar{\rho}_0:=\rho_Ay^A$, $\bar{\rho}_{10}:=\rho_{1A}y^A$,
 $\bar{\rho}_{00}:=\rho_{AB}y^Ay^B$, $\bar{\tau}_0:=\tau_Ay^A$.

 Under the local coordinate system $(s,y^A)$, (\ref{y25}) can be
 written in the form $A\bar{\alpha}^2+B\bar{\alpha}+C=0$. So we
 have $A\bar{\alpha}^2+C=0$ and $B=0$. By $A\bar{\alpha}^2+C=0$ we
 have
  \beq
 &&\bar{\rho}_0^2-\bar{\rho}_{00}=\delta\bar{\alpha}^2,\label{w28}\\
 &&3b^2\tau^2s^2(1+ks^2)^2(\phi-s\phi')^2+bs^2(1+ks^2)\big[2(m-1)b\tau^2-\tau_1\big]\phi(\phi-s\phi')\nonumber\\
 &&\ \ \ \ -Kb^2\phi^4+\Big\{\big[(m+1)b\tau_1-(m+1)^2b^2\tau^2+\mu-\delta\big]t^2+b^2\delta\Big\}\phi^2=0,\label{w30}
  \eeq
where $\delta=\delta(x)$ and $\mu=\mu(x):= \rho_1^2-\rho_{11}$ are
some scalar functions.

 We  consider (\ref{w30}). Plug the
Taylor expansion of (\ref{y5}) into (\ref{w30}) and let $p_i$ be
the coefficients of $s^i$ in (\ref{w30}).

\

\noindent {\bf Case I:} Assume $a_1\ne 0$. We will show this is
impossible.

 Firstly we have
$\delta=0$ from $p_2=0$. Then $p_4=0$ gives
 \be\label{w35}
 \mu=(m+1)^2b^2\tau^2-(m+1)b\tau_1+b^2a_1^2K.
 \ee
Plug $\delta=0$ and (\ref{w35}) into $p_{m+3}=0$ and we obtain
 \be\label{w36}
  \tau_1=2(m-1)b\tau^2+\frac{2ba_1^2K}{m-1}.
 \ee
Plug (\ref{y5}), $\delta=0$, (\ref{w35}) and (\ref{w36}) into
(\ref{w30}) and then it is clear that $p_{2m+2}=0$ gives
 \be\label{w37}
 K=\frac{(m-1)^2}{a_1^2}\tau^2,
 \ee
and $p_{3m+1}=0$ gives $K=0$. Thus we have $\tau=0$ from
(\ref{w37}). Then by (\ref{cw6}) we get a contradiction.

\

\noindent {\bf Case II:} As shown above we have $a_1= 0$. Then
$p_{2m}=0$ gives $\delta=0$. Plug $\delta=0$ into $p_{2m+2}=0$ and
we have
 \be\label{w38}
  \mu=4mb^2\tau^2-2mb\tau_1.
 \ee
Now plugging $\delta=0$ and (\ref{w38}) into (\ref{w30}) yields
$K=0$.     \qed

\section{The third class in Theorem \ref{th4}}
In this section, we study the property of an
$(\alpha,\beta)$-metric determined by (\ref{y16}), (\ref{w06}) and
(\ref{ycr06}) in Theorem \ref{th4}.

\begin{prop}
 Let $F$ be an $n$-dimensional ($\alpha,\beta$)-metric given by the $\phi(s)$ in (\ref{y16}) which is  projectively flat with
 constant flag curvature $K$.
  Then we have $K=0$ and $F$ is locally Minkowskian.
\end{prop}

{\bf Proof.} The projective factor $P$ is given by (\ref{w03}). By
(\ref{w06}) we get
 \be\label{y0022}
 r_{00}=2\tau\big(mb^2-(m+1+kb^2)s^2\big)\alpha^2-\frac{2(m+1+2kb^2)}{(m-1)b^2}s\alpha s_0.
 \ee
By (\ref{ycr06}) we get
 \be\label{y0023}
 \alpha_0:=\alpha_{x^k}y^k=a_{ir}y^r\frac{2}{\alpha}G^i_{\alpha}
 =2\rho\alpha+\frac{2(k\beta^2-m\alpha^2)\big[(m-1)\tau
 b^2\beta+s_0\big]}{(m-1)b^2\alpha}.
 \ee
Now it follows from (\ref{ycr06}), (\ref{y0022}) and (\ref{y0023})
that
 \be\label{y0024}
 s_{x^k}y^k=-\frac{s\alpha_0}{\alpha}+\frac{r_{00}+2b_iG^i_{\alpha}}{\alpha}
           =-\frac{2\beta(\alpha^2+k\beta^2)\big[(m-1)\tau b^2\beta+s_0\big]}{(m-1)b^2\alpha^3}.
 \ee
  Then plug (\ref{y16}), (\ref{w03}), (\ref{y0023}) and
(\ref{y0024}) into (\ref{y4}) and we obtain
 \be
(m-1)^2b^4K\beta^{2m}(\alpha^2+k\beta^2)^{1-m}+A_0=0,\label{y45}
 \ee
where $A_0$ is a polynomial in $(y^i)$. Since $m\ne 0,1$, clearly
we get $K=0$ from (\ref{y45}). Further, we conclude that $F$ is
Berwaldian since by (\ref{w03}), the projective factor $P$ is a
1-form. Thus $F$ is locally Minkowskian.        \qed

\section{The fourth class in Theorem \ref{th4}}
In this section, we study the property of the
$(\alpha,\beta)$-metric determined by (\ref{j4}), (\ref{j5}) and
(\ref{j6}) in Theorem \ref{th4}.

\begin{prop}\label{p51}
Let $F=\alpha \phi(s)$, $s=\beta/\alpha$, be an
  $(\alpha,\beta)$-metric on an open subset $U$ of the two-dimensional
  Euclidean space $R^2$. Suppose that
    $\beta$ is not parallel with respect to $\alpha$ and $F$ satisfies (\ref{j4}), (\ref{j5}) and
(\ref{j6}) with constant flag curvature $K$. Then we have one of
the following cases:
 \ben
 \item[{\rm (i)}] $\phi$ is given by
   \be\label{w088}
  \phi(s)=s^m\big(1-\frac{s^2}{b^2}\big)^{\frac{1-m}{2}}.
  \ee
  In this case, we
   have $K=0$.

\item[{\rm (ii)}] $\phi$ is given by
 \be\label{y088}
 \phi(s)= \frac{2b}{1-kb^2}\Big\{b\sqrt{1-ks^2}
         -\sqrt{b^2-s^2}\Big\},
 \ee
 where $k\ne 1/b^2$.

\item[{\rm (iii)}] $\phi$ is given by
 \be\label{y089}
  \phi(s)=\frac{4b^2}{(1-kb^2)^2}\frac{\big(b\sqrt{1-ks^2}
          -\sqrt{b^2-s^2}\big)^2}{\sqrt{1-ks^2}},
 \ee
where $k\ne 1/b^2$. In this case, we  have $K=0$.
 \een
\end{prop}

{\bf Proof.} The projective factor $P$ is given by (\ref{w04}). By
(\ref{j5}) and (\ref{j6}) we get
 \be\label{w50}
  r_{00}=-\frac{2s}{b^2}\alpha s_0, \ \ \alpha_{x^k}y^k=
  \Big\{2\rho-\frac{2(m-2+kt^2)s_0}{(m-1)b^2}\Big\}\alpha, \ \
  s_{x^k}y^k=\frac{2s(ks^2-1)}{(m-1)b^2}s_0.
 \ee
It follows from (\ref{j4}) that $\phi(s)$ satisfies
 \be\label{w51}
 \frac{\phi-s\phi'+(b^2-s^2)\phi''}{s\phi+(b^2-s^2)\phi'}=\frac{m-1}{s(1-ks^2)}.
 \ee
Now substitute (\ref{w50}) and (\ref{w51}) into (\ref{y4}) and we
obtain
 \beq
 K\alpha^2\phi^2\hspace{-0.6cm}&&=\rho^2-\rho_{00}+\frac{1}{b^2}S_{00}
 -\frac{2\rho}{b^2}s_0
 +\frac{3(ks^2-1)^2}{(m-1)^2b^4\phi^2}(\phi-s\phi')^2s_0^2\nonumber\\
 &&+\frac{ks^2-1}{(m-1)b^2\phi}\Big\{S_{00}
 +\frac{2(kb^2s^2-2s^2+b^2)}{b^2(b^2-s^2)}s_0^2
 -2\rho s_0\Big\}(\phi-s\phi')\nonumber\\
&&+\frac{(1-m+2kb^2)s^2+(m-3)b^2}{(m-1)(b^2-s^2)b^4}s_0^2,\label{y67}
 \eeq
where $\rho_{00}$ and $S_{00}$ are defined by
 $$
\rho_{ij}:=\frac{1}{2}(\frac{\pa \rho_i}{\pa x^j}+\frac{\pa
  \rho_j}{\pa x^i}),\ \ \rho_{00}:=\rho_{ij}y^iy^j,\ \ S_{ij}:=\frac{1}{2}(\frac{\pa s_i}{\pa x^j}+\frac{\pa
  s_j}{\pa x^i}),\ \ S_{00}:=S_{ij}y^iy^j.
 $$
 Since the
dependence of $\phi$ on $s$ is not clear, we choose a special
coordinate system $(s,y^2)$ at a fixed point $x_o$ as that in
Section \ref{sec3}.

Under the local coordinate system $(s,y^2)$, put $\rho$ and
$\rho_{00}$ as in (\ref{y27}) and

 $$
 S_{00}=\frac{S_{11}s^2}{b^2-s^2}\bar{\alpha}^2
  +\frac{2s\bar{S}_{10}}{\sqrt{b^2-s^2}}\bar{\alpha}+\bar{S}_{00}.
 $$
Substitute them into (\ref{y67}) and then (\ref{y67}) can be
written in the form
  \be\label{y68}
 A\bar{\alpha}^2+B\bar{\alpha}+C=0,
  \ee
 where $A,B,C$ are some polynomials in $(y^2)$.
By (\ref{y68}) we have $A\bar{\alpha}^2+C=0$ and $B=0$. Since $F$
is two-dimensional, $A\bar{\alpha}^2+C=0$
 can be written as
 \beq
 &&3(b^2-s^2)(1-ks^2)^2(\phi-s\phi')^2s_2^2-(1-ks^2)\Big\{\big[(2kb^2+b^2-mb^2-4)s_2^2\nonumber\\
 &&\ \ +(m-1)b^2(2\rho_2s_2+\eta+S_{11})\big]s^2+b^2(mb^2-b^2+2)s_2^2\nonumber\\
 &&-(m-1)b^4(2\rho_2s_2+\eta)\Big\}\phi(\phi-s\phi')+(m-1)\phi^2\Big\{\big[(2kb^2-mb^2+b^2-m+1)s_2^2\nonumber\\
 &&\ \
 +(m-1)b^2(2\rho_2s_2+\eta+S_{11}+b^2\mu-b^2\delta)\big]s^2+b^2(mb^2-b^2+m-3)s_2^2\nonumber\\
 &&\ \
 -(m-1)b^4(2\rho_2s_2+Kb^2\phi^2+\eta-b^2\delta)\Big\}=0,\label{y69}
 \eeq
where
 $$\mu:=\rho_1^2-\rho_{11},\ \ \delta:=\rho_2^2-\rho_{22},\ \
 \eta=s_2^2-S_{22},$$
Put
 \be\label{w055}
 \phi(s)=s^m+a_{m+2}s^{m+2}+a_{m+4}s^{m+4}+a_{m+6}s^{m+6}+a_{m+8}s^{m+8}+o(s^{m+8}),
 \ee
and plug it into (\ref{w51}). Then we obtain
 \beq
 a_{m+2}&=&\frac{(m-1)(mkb^2+2)}{2(m+2)b^2},\ \
 a_{m+4}=\frac{(m^2-1)[mkb^2(mkb^2+2kb^2+4)+8]}{8(m+2)(m+4)b^4},\label{w56}\\
 a_{m+6}&=&\frac{m+3}{6(m+6)b^2}\big\{(mkb^2+4kb^2+6)a_{m+4}-(m+1)ka_{m+2}\big\},\label{w57}\\
 a_{m+8}&=&\frac{m+5}{8(m+8)b^2}\big\{(mkb^2+6kb^2+8)a_{m+6}-(m+3)ka_{m+4}\big\}.\label{w58}
 \eeq
Substitute (\ref{w055}) into (\ref{y69}) and let $p_i$ be the
coefficients of $s^i$ in (\ref{y69}). Then we can easily get
$p_{2m}$ and $p_{2m+2}$ ($m\ge 2$), $p_{2m+4}$ ($m\ge 3$),
$p_{2m+6}$ ($m\ge 4$) and $p_{2m+8}$ ($m\ge 5$). Here we omit
their expressions.

\bigskip

\noindent {\bf Case I:} Assume $k=1/b^2$. Plug $k=1/b^2$ into
(\ref{w51}) and then solving it gives
 \be\label{w59}
 \phi(s)=s^m(1-\frac{s^2}{b^2})^{\frac{1-m}{2}}.
 \ee
Plug $k=1/b^2$ into $p_{2m}=0$ and $p_{2m+2}=0$ and we have
 \be\label{w60}
 \mu=-\frac{2}{b^2}S_{11},\ \
 \delta=\frac{2}{b^2}\big[b^2(\eta+2\rho_2s_2)-(2+b^2)s_2^2\big].
 \ee
Substitute $k=1/b^2$, (\ref{w59}) and (\ref{w60}) into (\ref{y69})
and we easily get $K=0$.

\bigskip

\noindent {\bf Case II:} Assume $k\ne 1/b^2$. We will show $m=2$
or $m=4$.

\bigskip

 {\bf Case II (A):} Suppose $m\ge 5$. This case will show
a contradiction.

Solve the system $p_{2m}=0$, $p_{2m+2}=0$, $p_{2m+4}=0$,
$p_{2m+6}=0$, and by using (\ref{w56}) and (\ref{w57}) we get

\be\label{w63}
 \mu=\frac{48(b^2k-1)}{(m+2)^2b^4}s_2^2,\ \ \ \
 \delta=\frac{(m-2)(5m^2+3m-32)}{(m^2-1)(m+2)b^4}s_2^2,
\ee
 \be\label{w64}
 S_{11}=\frac{(1-kb^2)(m^3+41m^2-10m+40)}{2(m^2-1)(m+2)^2b^2}s_2^2,
\ee
 \be\label{w65}
 \eta=\Big\{1+\frac{9m^3+m^2-42m+56}{2(m^2-1)(m+2)b^2}\Big\}s_2^2-2\rho_2s_2.
 \ee
Substitute (\ref{w63})--(\ref{w65}) into $p_{2m+8}=0$ and by using
(\ref{w56})--(\ref{w58}) we obtain $s_2=0$, which contradicts with
that $\beta$ is not parallel.

\bigskip

 {\bf Case II (B):} Put $m=3$. This case is also a
contradiction.

Likely, we can first get $\mu,\delta,\eta$ and $S_{11}$ by solving
the system
$$p_6=0,p_8=0,p_{10}=0,p_{12}=0.$$
Then plugging them into $p_{14}=0$ yields
 \be\label{w66}
 s_2^2=\frac{6125}{8(kb^2-1)^3}Kb^{10}.
 \ee
Finally, we obtain $K=0$ by $p_{16}=0$, where we need
(\ref{w56})--(\ref{w58}) ($m=3$) and
 $$
 a_{m+10}=\frac{m+7}{10(m+10)b^2}\big\{(mkb^2+8kb^2+10)a_{m+8}-(m+5)ka_{m+6}\big\}, \ \ (m=3).
 $$
Thus we get a contradiction by (\ref{w66}) and $K=0$.

\bigskip

 {\bf Case II (C):} Put $m=2$. This case gives
Proposition \ref{p51}(ii) since the ODE (\ref{w51}) with $m=2$
gives (\ref{y088}).

\bigskip

 {\bf Case II (D):} Put $m=4$. This case gives
Proposition \ref{p51}(iii) since the ODE (\ref{w51}) with $m=4$
gives (\ref{y089}).

\section{The fifth class in Theorem \ref{th4}}
In this section, we study the property of the
$(\alpha,\beta)$-metric determined by (\ref{yg5}), (\ref{yg6}) and
(\ref{cw002}) in Theorem \ref{th4}. Note that in Theorem
\ref{th4}(v), if $k_1-k_2^2=0$, then Theorem \ref{th4}(v) is a
special case of Theorem \ref{th4}(iii) with $m=-3$. So we only
need to assume $k_1-k_2^2\ne 0$ in the following proposition.

\begin{prop}
 Let $F=k_1\beta+2k_2\alpha^2/\beta+\alpha^4/\beta^3$ with $k_1-k_2^2\ne 0$ be
  a two-dimensional ($\alpha,\beta$)-metric  which is  projectively flat with
 constant flag curvature $K$.  Then we have
 $K=0$ and $F$ is locally Minkowskian. In this case,  $\alpha$ is flat and $\beta$ is parallel with respect
 to $\alpha$.
\end{prop}

{\bf Proof.}
 The projective factor $P$ is given by (\ref{cw0002}).
Similarly as the discussions in the above classes,  by (\ref{yg6})
we can get the expressions of $r_{00}$ and $r_0$, and by
(\ref{cw002}) we can get the expressions of
 $$
 \alpha_0:=\alpha_{x^k}y^k=a_{ir}y^r\frac{2}{\alpha}G^i_{\alpha},
 \ \ s_{x^k}y^k=-\frac{s\alpha_0}{\alpha}+\frac{r_{00}+2b_iG^i_{\alpha}}{\alpha}
 $$
Then plug them together with (\ref{yg5}) and (\ref{cw0002})
 into (\ref{y4}) and thus (\ref{y4}) can be firstly written in the
 following form
  \be\label{cr8}
 A_1\beta^2-64Kb^4(1+k_2b^2)^3\alpha^{16}=0,
  \ee
where $A_1$ is a homogeneous polynomial in $(y^i)$. Clearly by
(\ref{cr8}) we get $K=0$. By $K=0$, (\ref{cr8}) has the following
equivalent form
 \be\label{cr9}
 A_2(\alpha^4+2k_2\alpha^2\beta^2+k_1\beta^4)+48(1+k_2b^2)\beta^4(\alpha^2+k_2\beta^2)^4T^2=0,
 \ee
where $A_2$ is a homogeneous polynomial in $(y^i)$ and $T$ is
defined by
 $$
 T:=b^2s_0\alpha^2+\big[8b^2\tau(1+k_2b^2)\beta-(2+k_2b^2)s_0\big]\beta^2.
 $$
Since $k_1\ne k_2^2$, it is easy to show that
$(\alpha^2+k_2\beta^2)^4$ and
$\alpha^4+2k_2\alpha^2\beta^2+k_1\beta^4$ have no common factor.
Thus $\alpha^4+2k_2\alpha^2\beta^2+k_1\beta^4$ must be divided by
$T$ and so we get $T=0$. Now by $T=0$ and the definition of $T$,
we easily show that
 $$s_0=0,\ \ \tau=0.$$
  Further, the second formula in (\ref{ygjcw}) always hold in
  two-dimensional case. So we have $s_{ij}=0$ by $s_i=0$. Thus
  $\beta$ is parallel with respect to $\alpha$.  \qed

\section{Proof of Theorem \ref{th01}}

Let $F=\beta^m\alpha^{1-m}$ be an $m$-Kropina metric ($m\ne 0,1$).
For the case $m\ne -1$, the conclusion has been proved under
weaker condition that $F$ is only assumed to be Douglasian if
$n=2$ (see \cite{Y4}) and $F$ is only assumed to be locally
projectively flat if $n\ge 3$ (see \cite{Y3}).  Therefore, we only
need to assume $m=-1$, that is, $F$ is a Kropina metric.

By the proof in Section \ref{first}, we have (\ref{cr5}) and
(\ref{cr7}) with $k=0$. Then $F$ is locally projectively flat with
constant flag curvature if and only if
 \beq
 r_{00}&=&\mu\alpha^2, \ \ s_{ij}=\frac{b_is_j-b_js_i}{b^2},\label{cr90}\\
G^i_{\alpha}&=&\rho
y^i-\frac{r_{00}}{2b^2}b^i-\frac{\alpha^2}{2b^2}s^i,\label{cr91}\\
 \rho&=&\frac{\mu\beta+s_0}{b^2}.\label{cr92}
 \eeq
Next we apply a useful deformation on $\alpha$ and $\beta$ to
obtain the local expressions of $\alpha$ and $\beta$ based on
(\ref{cr90})--(\ref{cr92}). Define
 \be\label{cr93}
 \widetilde{\alpha}:=\frac{1}{b}\alpha,\ \
 \widetilde{\beta}:=\frac{1}{b^2}\beta.
 \ee
 In \cite{Y3} and \cite{Y4} we define a deformation applied to $m$-Kropina
 metric as follows
  $$
\widetilde{\alpha}:=b^m\alpha,\ \
 \widetilde{\beta}:=b^{m-1}\beta.
  $$
When $m=-1$, the above becomes (\ref{cr93}). Under the deformation
(\ref{cr93}), by (\ref{cr90})--(\ref{cr92}), we obtain
 \be\label{cr94}
 \widetilde{r}_{ij}=0,\ \ \widetilde{s}_{ij}=0,\ \ \widetilde{G}^i=0.
 \ee
Therefore, by (\ref{cr94}) we conclude that $\widetilde{\alpha}$
is flat and $\widetilde{\beta}$ is parallel with respect to
$\widetilde{\alpha}$. Now put $\eta:=b^2$, and then we get
(\ref{ycw2}) with $m=-1$.    \qed

\bigskip

{\bf Acknowledgement:}

The  author expresses his sincere thanks to China Scholarship
Council for its funding support.  He did this job  during the
period (June 2012--June 2013) when he as a postdoctoral researcher
visited Indiana University-Purdue University Indianapolis, USA.

\vspace{0.6cm}

\noindent Guojun Yang \\
Department of Mathematics \\
Sichuan University \\
Chengdu 610064, P. R. China \\
 ygjsl2000@yahoo.com.cn

\end{document}